\documentclass[12pt,leqno]{article}
\def\s{\dot{s}}

\def\int{\mathbb{Z}}

\def\Ue{{\cal U}_{\varepsilon}({\mathfrak g})}

\def\proof{{\bf Proof. }}

\def\pf{\proof}

\usepackage{times}
\usepackage{graphics}
\usepackage{amssymb}
\usepackage{amscd}
\usepackage{amsmath}
\input xy  
\xyoption{all}
\title{On Lusztig's map for spherical  unipotent conjugacy classes}  
\newtheorem{theorem}{Theorem}[section]
\newtheorem{lemma}[theorem]{Lemma}
\newtheorem{corollary}[theorem]{Corollary}
\newtheorem{proposition}[theorem]{Proposition}

\newtheorem{remark}[theorem]{Remark}

\author{Giovanna Carnovale, Mauro Costantini\\
Dipartimento di Matematica\\
via Trieste 63 - 35121 Padova - Italy\\
email: carnoval@math.unipd.it, costantini@math.unipd.it}

\date{}

\begin{document}
\maketitle
\begin{abstract}
We provide an alternative description of the restriction to spherical unipotent conjugacy classes, of Lusztig's map $\Psi$ from the set of unipotent conjugacy classes in a connected reductive algebraic group to the set of conjugacy classes of its Weyl group. For irreducible root systems, we analyze the image of this restricted map and we prove that a conjugacy class in a finite Weyl group has a unique maximal length element if and only if it has a maximum.
\end{abstract}


\section{Introduction}

In \cite{S}, Springer has shown how to associate to a unipotent conjugacy class of a connected reductive algebraic group $G$ over an algebraically closed field $k$ some irreducible representations of the associated Weyl group $W$. This lead Kazhdan and Lusztig to the definition, in \cite{KL}, of a conjecturally injective map from the set $\underline{G}$ of unipotent conjugacy classes of $G$ to the set $\underline{W}$ of conjugacy classes of $W$, for $k={\mathbb C}$. This map is not easily computable but Lusztig has very recently introduced in \cite{lusztig,lusztig2} a new, more computable, surjective map $\phi$ defined in all characteristics, from $\underline{W}$ to $\underline{G}$, and a right inverse $\Psi$  which conjecturally coincides with the Kazhdan-Lusztig map over the complex numbers.  The map $\phi$ is defined by assigning to a conjugacy class $C$ in ${W}$ a minimal unipotent conjugacy class in $G$, with respect to Zariski closure, having non-empty intersection with the Bruhat double coset corresponding to a minimal length element in $C$. It is a non-trivial result that this construction works.  The proof of this important property  is split into a proof  for classical groups and one, based on a computer calculation, for exceptional ones. The right inverse $\Psi$ is defined by taking, for a given unipotent class $\gamma$ in $G$, the unique class $C$ in $W$ in the fiber of $\gamma$ through $\phi$ for which the dimension of the fixed point space of $w\in C$ in the geometric representation of $W$  is minimal. Also in this case, the fact  that this procedure actually works is a deep result. 

The aim of this note is to give a different and direct combinatorial description of the restriction to spherical unipotent conjugacy classes of the map $\Psi$. 
We recall that a conjugacy class $C$ in $G$ is called spherical if a Borel subgroup $B$ of $G$ has a dense orbit in $C$.
This new description is made possible by several recent results showing how the  relation between spherical conjugacy classes and the Bruhat decomposition can be made very explicit. It has been shown in \cite{ccc,gio-mathZ,mauro-cattiva,lu} that spherical (unipotent) classes may be characterized by means of a dimension formula involving the maximal Weyl group element $w$ for which $BwB$ meets a class.
More precisely, let us  define, for  $\gamma$ in $\underline{G}$, the element $w_\gamma\in W$ as the unique element in $W$ for which $Bw_\gamma B\cap \gamma$ is Zariski dense in $\gamma$. Then, $\gamma$ is spherical if and only if  $\dim \gamma=\ell(w_\gamma)+{\rm rk}(1-w_\gamma)$, where $\ell$ is the length function on $W$ and ${\rm rk}$ is the rank of the operator in the geometric representation of $W$. In addition, spherical conjugacy classes in good, odd characteristic are also characterized as those classes intersecting only Bruhat double cosets corresponding to involutions (\cite{gio-mathZ,gio-fourier}). Combining all these properties with the analysis of the elements $w_\gamma$ in \cite{chan-lu-to} leads us to the proof of our main result:

\medskip

\noindent{\bf Theorem}\;{\em Let $\gamma$ be a spherical unipotent conjugacy class. Then, $\Psi(\gamma)=W\cdot w_\gamma$.}

\medskip

We also give some results on the map $\iota\colon{\underline{G}}\to\underline{W}$ defined by $\iota(\gamma)=W\cdot w_\gamma$. 

This map can be defined on the set of all conjugacy classes in $G$. It was observed in \cite[Remark 3]{chan-lu-to} that the image of  the set of all conjugacy classes and of the set of all spherical conjugacy classes through this map, in characteristic zero or good and odd characteristic,  is the set $\underline{W}_m$ of classes in $\underline{W}$ having a unique maximal length element. 
We analyze the image of the restriction of $\iota$ to the set $\underline{G}_{sph}$ of spherical unipotent conjugacy classes. A case-by-case analysis allows us to conclude  that

\medskip

\noindent{\bf Proposition}\;{\em For every irreducible root system there always exists a $p$ such that in characteristic $p$ we have $\iota(\underline{G}_{sph})=\underline{W}_m$. }

\medskip

It is worthwhile to mention that the element $w_\gamma$, for spherical classes, controls the $G$-module structure of the ring of regular functions ${\mathbb C}[\gamma]$. Indeed, this module is multiplicity-free by \cite{VK} and  it has been observed in \cite{ccc} that the weights $\lambda$ occurring in the decomposition of ${\mathbb C}[\gamma]$ all satisfy the equality $w_\gamma\lambda=-\lambda$ and that the rank of the lattice generated by these weights is ${\rm rk}(1-w_\gamma)$. The precise $G$-module decomposition of  ${\mathbb C}[\gamma]$ has been given in \cite{mauro-mathZ}.  

\medskip

We conclude the paper by proving that 

\medskip

\noindent{\bf Theorem}\;{\em The set $\underline{W}_m$  coincides with the set of classes in $\underline{W}$ having  maximum element with respect to the Bruhat order.}

\medskip

This result holds for arbitrary finite Coxeter groups (see Remark~\ref{hultman}).
 
\section{Notation}

Throughout this paper $G$ is a semisimple algebraic group over an algebraically closed field $k$. 
Let $T$ be a maximal torus of $G$,  and let $\Phi$ be the associated
root system. Let $B\supset T$ be a Borel subgroup, $B^-$ its opposite Borel subgroup, and let
$\Delta=\{\alpha_1,\ldots,\alpha_n\}$ be the basis of $\Phi$ relative
to $(T,\,B)$.  
The Weyl group is  denoted by $W=N(T)/T$,
the symbol $\underline{W}$ will indicate the set of its conjugacy classes, and $\underline{W}_{inv}$ will indicate the set of conjugacy classes of involutions in $W$, that is, the set of classes of those elements $w\in W$ such that $w^2=1$.
The symbol $\underline{G}$ will stand for the set of unipotent conjugacy classes and $\underline{G}_{sph}$ will denote the set of spherical unipotent ones. We recall that a conjugacy class $\gamma$  in $G$ is called spherical if $B$ has a Zariski dense orbit in $\gamma$. 

For any $C\in\underline{W}$ we define $C_{min}$ to be the subset of $C$ consisting of  elements of minimal length.
For $w\in W$ we define $\Sigma_w=\{\gamma\in\underline{G}~|~\gamma\cap BwB\neq\emptyset\}$. 

For $\gamma\in\underline{G}$ we define $W_\gamma=\{w\in W~|~\gamma\cap BwB\neq\emptyset\}$. It is clear that $W_\gamma$ is always not empty. It is also true that $\Sigma_w$ is always not empty: indeed $BwB\cap B^-\neq\emptyset$ for every $w\in W$ (\cite[\S A2]{HL}), so $BwB\cap U^-\neq\emptyset$ for every $w\in W$.

As usual, $w_0$ denotes the longest element in $W$ and, for $\Sigma\subseteq\Delta$, we shall denote by $w_\Sigma$ the longest element in the parabolic subgroup $W_\Sigma$ of $W$ generated by simple reflections indexed by elements in $\Sigma$. The root subsystem of $\Phi$ generated by the roots in $\Sigma$ will be denoted by $\Phi_\Sigma$.

It follows from \cite[8.2.6(b)]{GP} and \cite[1.2(a)]{lusztig}
that for $w,\sigma\in C_{min}$ then $\Sigma_w=\Sigma_\sigma$. 

Let $\phi\colon\underline{W}\to \underline{G}$ be the map introduced in \cite{lusztig}. It is defined as follows: let $C\in\underline{W}$ and let $w\in C_{min}$. The image of $C$ through $\phi$ is the unique $\gamma\in\underline{G}$ such that $\gamma\in\Sigma_w$ and such that every $\gamma'\in\underline{G}$  lying in $\Sigma_w$ contains $\gamma$ in its closure. By \cite[Theorem 0.4]{lusztig} the map $\phi$ is surjective. 

If $\gamma\in\underline{G}$  and $C\in \phi^{-1}(\gamma)$ then $\gamma\in \Sigma_w$ for some $w\in C_{min}$. For $\gamma$ a spherical unipotent conjugacy class the set $W_\gamma$ has a particular structure. We recall the facts we will need.

\begin{theorem}\label{p=2}(\cite{gio-mathZ,mauro-cattiva}) Let $\gamma$ be a spherical conjugacy class, and let $\gamma\cap BwB$ be non-empty. Assume in addition that $\gamma$ is unipotent if ${\rm char}(k)=2$. Then, $w$ is an involution. 
\end{theorem}
\pf  If ${\rm char}(k)$ is zero or good and odd this is \cite[Theorem 2.7]{gio-mathZ}. The  same proof holds as long as ${\rm char}(k)\neq2$.  For ${\rm char}(k)=2$, let $u$ be an element of $\gamma\cap BwB$. From the classification of spherical unipotent conjugacy classes it follows that $u$ is an involution, see \cite[Theorem 3.18]{mauro-cattiva}. Thus, $u=u^{-1}\in Bw^{-1}B\cap BwB$, forcing $w=w^{-1}$.\hfill$\Box$ 

\bigskip

So, $\phi^{-1}(\underline{G}_{sph})\subseteq \underline{W}_{inv}$. One may wish to see whether  $\underline{G}_{sph}$ can be characterized as the image of a suitable subset of $\underline{W}_{inv}$.

\medskip

The statement of  the Lemma below was communicated to the first named author by Kei-Yuen Chan.
\begin{lemma}\label{chan}Let ${\rm char}(k)\neq2$. Let  $\gamma$ be a (not necessarily unipotent) spherical conjugacy class and let $\gamma\cap BwB\neq\emptyset$ for some $w\in C$ and $C\in\underline{W}$. Then $\gamma\cap B\sigma B\neq\emptyset$ for every $\sigma\in C$. The same conclusion holds for ${\rm char}(k)=2$ if $\gamma$ is a spherical unipotent conjugacy class. 
\end{lemma}
\pf Let $\sigma=s_{i_l}\cdots s_{i_1}w s_{i_1}\cdots s_{i_l}$  with $\tau=s_{i_l}\cdots s_{i_1}$ of minimal length $l$ such that $\sigma=\tau w\tau^{-1}$. Let us put
$\sigma_j=s_{i_j}\cdots s_{i_1}w s_{i_1}\cdots s_{i_j}$ for $j=0,\ldots l$, so that $\sigma_0=w$ and $\sigma_l=\sigma$. We shall prove by induction on $j$ that $\gamma\cap B\sigma_j B\neq\emptyset$ for every $j\in\{0,\,\ldots,\,l\}$. The basis of the induction is our assumption. Assume $\gamma\cap B\sigma_j B\neq\emptyset$ for a given $j$. Then, there is also
 $x\in B\sigma_j\cap \gamma$.  Let $\s_{i_{j+1}}$ be a lift of $s_{i_{j+1}}$ in $N(T)$. We have 
$$\s_{i_{j+1}}x\s_{i_{j+1}}^{-1}\in s_{i_{j+1}}B\sigma_js_{i_{j+1}}\subseteq B\sigma_{j+1}B\cup B \sigma_j s_{i_{j+1}}B.$$ 
By Theorem~\ref{p=2} the class $\gamma$ intersects only cells corresponding to involutions. Hence, $w$ and $\sigma_j$ are involutions. On the other hand, $\sigma_j s_{i_{j+1}}$ is an involution if and only if $\sigma_j$ and $s_{i_{j+1}}$ commute, but this would contradict minimality of the length of  $\tau$. Thus,
$\gamma\cap B \sigma_j s_{i_{j+1}}B=\emptyset$, and we necessarily have $\s_{i_{j+1}}x\s_{i_{j+1}}^{-1}\in B\sigma_{j+1}B\cap \gamma$, yielding the statement.\hfill$\Box$

\bigskip

Let $\gamma$ be any conjugacy class in $G$. We shall  denote by $w_\gamma$ the unique element in $W$ for which $B w_\gamma B\cap \gamma$ is dense in $\gamma$, and by $C^\gamma=W\cdot w_\gamma$, the conjugacy class of $w_\gamma$ in $W$.  
Let us denote by $\underline{W}_m$ the set of classes in $\underline{W}$ containing a unique maximal length element.
We recall some basic facts.

\begin{theorem}\label{clt}(\cite{chan-lu-to}). Let $\gamma$ be a conjugacy class in $G$ and let $w_\gamma$ and $C^\gamma$ be as above. Then
\begin{enumerate}
\item $C^\gamma$ lies in $\underline{W}_m$ and $w_\gamma$ is its maximal length element;
\item $\underline{W}_m\subseteq\underline{W}_{inv}$; 
\item If ${\rm char}(k)$ is either $0$ or good and odd, then for every $C\in\underline{W}_m$ there is a spherical conjugacy class $\gamma$ such that 
$C=C^\gamma$. 
\end{enumerate}
\end{theorem}
\pf Statement 1. is  Corollary 2.11 in {\it loc. cit.}, the proof of which is characteristic-free. Statement 2. follows from the fact that any $w$ is conjugate to $w^{-1}$ (\cite[Theorem C]{carter-cc}). Statement 3. is observed in Remark 3 in {\it loc. cit.} 
\hfill$\Box$

\bigskip

We will also make use of the following result
\begin{theorem}\label{dimension}(\cite{ccc,gio-mathZ,mauro-cattiva,lu}) Let $\gamma$ be a unipotent conjugacy class, let $w_\gamma$ be as above, and let $w\in W$. 
\begin{enumerate}
\item If $\gamma\in\Sigma_w$ then $\dim \gamma\geq\ell(w)+{\rm rk}(1-w)$;
\item $\dim \gamma\geq\ell(w_\gamma)+{\rm rk}(1-w_\gamma)$;
\item $\gamma$ is spherical if and only if $\dim \gamma=\ell(w_\gamma)+{\rm rk}(1-w_\gamma)$.
\end{enumerate}
\end{theorem} 
 
\begin{proposition}\label{phi}
Let $\gamma$ be a spherical unipotent conjugacy class and let $C^\gamma$ be as above. Then $\phi(C^\gamma)=\gamma$.
\end{proposition}
\pf  
Let $w\in (C^\gamma)_{min}$. We need to show that $\gamma\in\Sigma_w$ and that it is the unique minimal element therein. 
 
By construction $\gamma$ lies in $\Sigma_{w_\gamma}$ so by Lemma~\ref{chan}, it also lies in $\Sigma_w$.  It follows from \cite[Propositions 2.8, 2.9]{chan-lu-to}, which in turn uses \cite[\S 2.9]{GKP} and  \cite[Proposition 3.4]{EG}, that if $\sigma\in C^\gamma$ and $y$ is a maximal length element in $C^\gamma$, then $\Sigma_\sigma\subseteq\Sigma_{y}$. In particular, this holds for $\sigma=w$ and $y=w_\gamma$ by Theorem~\ref{clt} (1). 

Let  $\gamma'\in\Sigma_w$. Then $\gamma'\in\Sigma_{w_\gamma}$ and by part 1 of Theorem~\ref{dimension}  we have $\dim \gamma'\geq \ell(w_\gamma)+{\rm rk}(1-w_\gamma)$. However, by Theorem~\ref{dimension} we have $\dim \gamma=\ell(w_\gamma)+{\rm rk}(1-w_\gamma)$
so $\gamma$ is minimal in $\Sigma_{w_\gamma}$,  and, {\it a fortiori}, in $\Sigma_w$. The assertion follows from  uniqueness of the minimal element in $\Sigma_w$ (see \cite{lusztig}).\hfill$\Box$

\bigskip

The above result can be rephrased by saying that the restriction to $\underline{G}_{sph}$ of the map 
$$\begin{array}{rl}
\iota\colon \underline{G}&\to \underline{W}_{inv}\\
\gamma&\mapsto C^\gamma
\end{array}$$
is a right inverse for $\phi$ on $\underline{G}_{sph}$. 

\bigskip

In \cite[Theorem 0.2]{lusztig2} a right inverse $\Psi$ to $\phi$ has been constructed.
It is defined as follows. For any $\gamma\in\underline{G}$ one considers $\phi^{-1}(\gamma)$. This set contains a unique element $C_0\in\underline{W}$ for which the dimension $d_C$ of the fixed point space of an (thus any) element in $C$ is minimal. Then $\Psi(\gamma)=C_0$. We want to compare the maps $\iota$ and $\Psi$ on $\underline{G}_{sph}$.

\bigskip

It is shown in \cite[Lemma 3.2]{chan-lu-to}  that $w_\gamma=w_0w_\Sigma$ for some $\Sigma\subseteq \Delta$ such that $w_\Sigma$ coincides with $w_0$ on $\Sigma$. Using the same arguments one can prove the following result, that we report here for completeness.

\begin{lemma}\label{maximal}Let $\gamma$ be a spherical unipotent conjugacy class or any spherical conjugacy class if ${\rm char}(k)$ is either $0$ or good and odd, and let $\sigma\in W_\gamma$ be a maximal length element in its conjugacy class $C$. Then,  $\sigma=w_0w_\Sigma$ for some $\Sigma\subseteq \Delta$ such that $w_\Sigma$ coincides with $w_0$ on $\Sigma$.
\end{lemma}
\pf Since $W_\gamma$ consists of involutions we may apply \cite[Theorem 1.1 (ii)]{perkins-rowley}, so $\sigma=w_0w_\Sigma$ for some $\Sigma\subseteq \Delta$.  In addition, $w_0$ and $w_\Sigma$ necessarily commute so $(-w_0)\Sigma=\Sigma$. Let $\alpha\in\Sigma$. We have $\beta=w_0w_\Sigma\alpha\in\Sigma\subseteq\Phi^+$ so $\ell(w_0w_\Sigma s_\alpha)=\ell(w_0w_\Sigma)+1$.
Then, by maximality of the length of $\sigma$ in $C$, we have $\ell(s_\alpha w_0w_\Sigma s_\alpha)=\ell(w_0w_\Sigma)$. By \cite[Lemma 3.2]{results}  we get $\alpha=\beta$. \hfill$\Box$

\bigskip

%
%
%

\begin{lemma}\label{rank}
Let $\Pi\subseteq \Delta$ and let $w=w_0w_\Pi$ be an involution with the property that $w_0$ restricted to $\Phi_\Pi$ is $w_\Pi$. Then, $(-w_0)(\Pi)=\Pi$ and
$${\rm rk}(1-w_0)={\rm rk}(1-w_\Pi)+{\rm rk}(1-w).$$
\end{lemma}
\pf The first statement follows from $w_0w_\Pi(\alpha)=\alpha$ for every $\alpha\in\Pi$. 

Let us denote by $E_{m}(x)$ the $m$-eigenspace of an operator $x$.
Clearly, if $x$ is an involution then $\dim E_{-1}(x)={\rm rk}(1-x)$.  It is an immediate exercise in linear algebra that if $x$ and $y$ are commuting involutions,  then $\dim E_{-1}(x)+\dim E_{-1}(y)=\dim E_{-1}(xy)$ if and only if $E_{-1}(x)\cap E_{-1}(y)=\{0\}$. 
%

We have $\Pi\subseteq E_1(w_0w_\Pi)=E_{-1}(w_0w_\Pi)^\perp$ so, since $w_\Pi$ can be written as a product of reflections with respect to roots in $\Pi$, for every $v\in E_{-1}(w_0w_\Pi)$ we have $w_\Pi(v)=v$. In other words,
$$E_{-1}(w_0w_\Pi)\cap E_{-1}(w_\Pi)\subseteq E_{1}(w_\Pi)\cap E_{-1}(w_\Pi)=\{0\}$$
whence the second statement.
%
%
%
%
 \hfill$\Box$

\bigskip


We are ready to state the main result of this paper.

\begin{theorem}\label{main}Lusztig's map $\Psi$  coincides with $\iota$ on $\underline{G}_{sph}$.
\end{theorem}
\pf Let $\gamma\in\underline{G}_{sph}$. By Proposition~\ref{phi} we have $C^\gamma\in\phi^{-1}(\gamma)$, so we only need to show that the dimension $d_C$ of the fixed point space $E_1(w)$ of an element $w\in C$ for $C\in \phi^{-1}(\gamma)$ is minimal for $w\in C^\gamma$.  

Let $C$ be a class in $\phi^{-1}(\gamma)$.  Then, every $\sigma$ in $C$ lies in $W_\gamma$ by Lemma~\ref{chan}.
%
By Theorem~\ref{p=2}, the set  $W_\gamma$ is a union of classes in $\underline{W}_{inv}$.  Moreover, all elements in $W_\gamma$ are less than or equal to $w_\gamma$ in the Bruhat ordering, in particular this holds for all elements in $C$.  
Let $z$ be a maximal length element in $C$. By Lemma~\ref{maximal},  $z=w_0w_\Sigma$ and $w_\gamma=w_0w_\Pi$ where $\Sigma$ and $\Pi$ are subsets of $\Delta$ on which $z$ and $w_\gamma$, respectively, act as the identity, and $z\leq w_\gamma$, or, equivalently,  $w_\Pi\leq w_\Sigma$. Since $w_\Sigma$ has a reduced expression as a product of reflections with respect to roots in $\Sigma$, the simple reflections occurring in some reduced expression of $w_\Pi$ correspond to some simple roots in $\Sigma$ by \cite[Corollary 2.2.3]{bb}.  By \cite[Corollary 1.4.8(ii)]{bb} the set of simple roots occurring in any reduced expression of $w_\Pi$  is precisely $\Pi$. Hence, $\Pi\subseteq \Sigma$. Moreover, the restriction of $w_\Sigma$ to $\Pi$ coincides with $w_\Pi$ so by Lemma \ref{rank} applied to $\Phi_\Sigma$ we have ${\rm rk}(1-w_\Sigma)={\rm rk}(1-w_\Pi w_\Sigma)+{\rm rk}(1-w_\Pi)$ so  ${\rm rk}(1-w_\Pi)\leq {\rm rk}(1-w_\Sigma)$. Applying Lemma \ref{rank} once more we see that 
$$\begin{array}{rl}
{\rm rk}(1-w_\gamma)&={\rm rk}(1-w_0w_\Pi)={\rm rk}(1-w_0)-{\rm rk}(1-w_\Pi)\\
&\geq {\rm rk}(1-w_0)-{\rm rk}(1-w_\Sigma)={\rm rk}(1-z)
\end{array}$$
so ${\rm rk}(1-x)$  reaches its maximum over $\phi^{-1}(\gamma)$ at $x=w_\gamma$. Since all elements in $\phi^{-1}(\gamma)$ are involutions, this gives precisely minimality of  $d_{C^\gamma}=\dim E_1(w_\gamma)$. Thus, $\Psi(\gamma)=C^\gamma$.\hfill$\Box$

%

\begin{corollary}The map  $\iota$ is injective on spherical unipotent conjugacy classes.
\end{corollary}

\begin{remark}{\rm Except for type $A_1$, the maps $\iota$ and $\Psi$ do not coincide on the full set $\underline{G}$ because $\Psi$ is necessarily injective whereas $\iota$ is not.  Indeed, the regular unipotent class $\gamma_{reg}$ intersects every $BwB$ (se \cite{Ka} or the result of Springer in the Appendix of \cite{EG}), so $\iota(\gamma_{reg})=W\cdot w_0$. On the other hand, there is always a spherical unipotent conjugacy class intersecting $Bw_0B$.}
\end{remark}

\bigskip

An important feature of the maps $\phi$ and $\Psi$ is that they are defined in all characteristic and they satisfy compatibility conditions as follows. 
For a fixed irreducible root system $\Phi$, let  $G^p$  denote a corresponding group in characteristic $p$ and let $\phi_p$, $\Psi_p$ and $\iota_p$ denote the corresponding maps $\phi$, $\Psi$ and $\iota$.  If in the sequel reference to $p$ is omitted, we shall mean that the statement holds for every $p\geq0$.  Let us recall that there is a dimension-preserving and order-preserving injective map $\pi\colon \underline{G^0}\to\underline{G^p}$ where the order is given by inclusion of Zariski closures (\cite[III, 5.2]{spaltenstein},\cite[\S 3.1]{lusztig2}). 
It is shown in \cite[Theorem 0.4(b)]{lusztig2} that $\Psi_0=\Psi_p\pi$ and $\pi=\phi_p\Psi_0$. The compatibility behaves well when we restrict to spherical conjugacy classes:

\begin{proposition}\label{pi}The map $\pi$ maps $\underline{G^0}_{sph}$ into $\underline{G^p}_{sph}$, and if $\gamma$ lies in $\underline{G^0}_{sph}$, then $w_{\pi(\gamma)}=w_\gamma$.\end{proposition}
\pf Let $\gamma\in\underline{G^0}_{sph}$. Then 
$$\pi(\gamma)=\phi_p\Psi_0(\gamma)=\phi_p\iota_0(\gamma)=\phi_p(C^\gamma).$$
Let $\sigma$ be a minimal length element in $C^\gamma$. Then,
$\pi(\gamma)\in \Sigma_\sigma$ and, arguing as in the proof of Lemma~\ref{phi}, since $w_\gamma$ is the maximal length element, $\pi(\gamma)\in\Sigma_{w_\gamma}$. Thus, $w_\gamma\leq w_{\pi(\gamma)}$. It is not hard to show, by induction on the length of a word in $W$, that if $w\leq \tau$ in the Bruhat order, then $\ell(w)+{\rm rk}(1-w)\leq\ell(\tau)+{\rm rk}(1-\tau)$ (see the proof of \cite[Proposition 6]{ccc}). Therefore, invoking  part 2 of
Theorem~\ref{dimension} for $\gamma$ we have 
$\dim(\pi(\gamma))=\dim(\gamma)=\ell(w_\gamma)+{\rm rk}(1-w_\gamma)\leq \ell(w_{\pi(\gamma)})+{\rm rk}(1-w_{\pi(\gamma)})$. 
Applying Theorem~\ref{dimension} to $\pi(\gamma)$, we have the first statement. The second one is immediate from $\Psi_0=\Psi_p\pi$ and Theorem~\ref{main}.\hfill$\Box$

\bigskip

In the remainder of the paper we analyze the image of the restriction of $\Psi$ to spherical unipotent conjugacy classes.

By part 1 of Theorem~\ref{clt}, the image of the restriction of $\iota$ to $\underline{G}_{sph}$ lies in $\underline{W}_m$. We observe that the map $\iota$ can be defined in the same way for any conjugacy class. 

Identifying a class in $\underline{W}_m$ with its unique maximal length element, we may endow $\underline{W}_m$ with a poset structure from the Bruhat order of $W$.  Inclusion of Zariski closures induces a poset structure on the set of conjugacy classes in $G$ and on $\underline{G}$.

We observe that if for some conjugacy classes $\gamma,\,\gamma'$ we have $\overline{\gamma}\subseteq\overline{\gamma'}$, then
$$\emptyset\neq B m_\gamma B\cap\gamma\subseteq \overline{B m_\gamma B\cap\gamma}=\overline{\gamma}\subseteq \overline{\gamma'}=\overline{Bm_{\gamma'}B\cap \gamma'}\subseteq \overline{B m_{\gamma'}B}$$ so 
$m_\gamma\leq m_{\gamma'}$ in the Bruhat order and $\iota$ is order-preserving.

By Theorem~\ref{clt}, in zero or good and odd characteristic the image of the set of all spherical classes through $\iota$ is exactly $\underline{W}_m$.  Let us analyze the situation for spherical unipotent conjugacy classes. 

\begin{proposition}\label{surjective}For every $\Phi$ there is some $p$ such that $\iota_p(\underline{G^p}_{sph})$ is $\underline{W}_m$. 
\end{proposition}
\pf  The list of the maximal length representatives for all elements in $\underline{W}_m$ is given in \cite{chan-lu-to} in terms of subdiagrams of the Dynkin diagram, and it can be deduced from \cite{perkins-rowley1}. In zero or good and odd characteristic we have $\iota_p(\underline{G^p}_{sph})=\underline{W}_m$ precisely in type $A_n$, $n\geq1$; $D_n$, $n\geq4$;  $E_6$; $E_7$; $E_8$ (see \cite[Table 3]{ccc}, \cite{gio-pacific},\cite[Lemma 3.5]{chan-lu-to}).  From Proposition~\ref{pi} it follows that in these cases we have $\iota_p(\underline{G^p}_{sph})=\underline{W}_m$ also when $p$ is a bad prime or $p=2$.


In type $C_n$ (and $B_n$), $n\geq 2$, in characteristic 2 there are $n+\left[\frac n2\right]$ non-trivial spherical unipotent conjugacy classes (see \cite[3.1.2]{mauro-cattiva}) and therefore we have $\iota_2(\underline{G^2}_{sph})=\underline{W}_m$.



In type $F_4$, for $p=3$ the poset of spherical unipotent conjugacy classes is the same as the corresponding poset in good characteristic, while for $p=2$ we have $\iota_2(\underline{G^2}_{sph})=\underline{W}_m$ (see \cite[Table 6, 7]{mauro-cattiva}).

In type $G_2$, for $p=2$ the poset of spherical unipotent conjugacy classes is the same as the corresponding poset in good characteristic, while for $p=3$ we have $\iota_3(\underline{G^3}_{sph})=\underline{W}_m$ (see \cite[Table 8, 9]{mauro-cattiva}).
\hfill$\Box$





\bigskip

\begin{corollary}The following are equivalent
\begin{enumerate}
\item $\iota_p(\underline{G^p}_{sph})=\underline{W}_m$ for every $p\geq0$; 
\item $\iota_0(\underline{G^0}_{sph})=\underline{W}_m$;
\item The restriction of $\pi$ to $\underline{G^0}_{sph}$ is an isomorphism onto $\underline{G^p}_{sph}$ for every $p\geq0$.
\end{enumerate}
\end{corollary}
\pf The equivalence of 1. and 2. is immediate from $\iota_0=\iota_p\pi$.  Let us assume 1. By bijectivity of $\iota_0$ and injectivity of $\iota_p$ we have 
$$|\underline{G^p}_{sph}|\leq|\underline{W}_m|=|\underline{G^0}_{sph}|$$ so injectivity of $\pi$ implies 3. Finally, Proposition \ref{surjective}  shows that 3. implies 1.\hfill$\Box$

\medskip
\begin{remark}{\rm Let ${\mathcal J}$ be the set of subsets of $\Delta$ such that 
$$\underline{W}_m=\{W\cdot w_0w_J~|~ J\in{\mathcal J}\}.$$  We can identify $\underline{W}_m$ with $\mathcal J$ and the partial order on $\underline{W}_m$ becomes reverse inclusion of subsets in ${\mathcal J}$. We observe that for $J,\,K\in{\mathcal J}$ both $J\cap K$ and $J\cup K$ are in ${\mathcal J}$ and therefore $\underline{W}_m$ is a lattice. It can be proved by inspection that for every $p$ the order-preserving map $\iota_p$, restricted to $\underline{G^p}_{sph}$ is a poset isomorphism onto its image and that $\underline{G^p}_{sph}$ is always a lattice.}
\end{remark}
\medskip

\begin{theorem}\label{max}The set $\underline{W}_m$ is the set of conjugacy classes in $W$ having maximum with respect to the Bruhat order.
\end{theorem}
\pf By Theorem~\ref{clt} (3) or by Proposition~\ref{surjective} if $C$ lies in $\underline{W}_m$, then  $C=C^\gamma$ for 
 some spherical conjugacy class in $G^p$ for some $p$. By Lemma \ref{chan}, every $w\in C$ lies in $W_\gamma$ so it must be less than or equal to $w_\gamma$ in the Bruhat ordering. Thus, the maximal length element  in $C$ is the sought maximum in $C$. Conversely, if $C$ has maximum $\sigma$ with respect to the Bruhat ordering then $\sigma$ has maximal length in $C$. Hence, $\sigma$ is the unique maximal length element in $C$ because for any $\tau\in C$ different from $\sigma$ we have $\ell(\tau)<\ell(\sigma)$.\hfill$\Box$

\bigskip

\begin{remark}\label{hultman}{\rm It was kindly suggested to us by A. Hultman that the statement of Theorem~\ref{max} for arbitrary finite Coxeter groups follows from  the observation in \cite[p. 577]{GKP}. Indeed, it is shown therein that for $C\in \underline{W}$ and any $w\in C$ there exists some $\sigma$ of  maximal length in $C$ and a chain of simple reflections $s_{i_1},\ldots,\,s_{i_r}$ satisfying 
$$\sigma_0=\sigma;\quad\sigma_{j}=s_{i_j}\sigma_{j-1}s_{i_j};\quad\sigma_r=w$$ and $\ell(\sigma_j)\geq\ell(\sigma_{j+1})$ for $j=0,\,\ldots,\,r$.
Now, if $C\in\underline{W}_m$ then $C\in \underline{W}_{inv}$ (see \cite[Thm. 8.7]{h3h4} for $H_3$ and $H_4$ or  \cite[Corollary 3.2.14]{GP} for arbitrary finite Coxeter groups). By \cite[Lemma 3.2]{results} we have $\ell(\sigma_j)=\ell(\sigma_{j+1})$ if and only if $\sigma_j=\sigma_{j+1}$, 
and  if $\ell(\sigma_j)>\ell(\sigma_{j+1})$, we necessarily have  $\ell(\sigma_j)=\ell(\sigma_{j+1})+2$. This forces
$$\sigma_{j-1}\geq \sigma_{j-1}s_{i_j}\geq s_{i_j}\sigma_{j-1}s_{i_j}=\sigma_j$$ in the Bruhat order, so the unique maximal length element $\sigma$ is the sought maximum in $C$. 

The main result in \cite{GKP} is based on a case-by-case analysis, but a new case-free proof is available in \cite{HN}. On the other hand,  surjectivity of $\iota$ on $\underline{W}_m$ relies on the case-by-case analysis in \cite{ccc}. This could be shortened by looking at the image through $\iota$ of the classes of involutions (in the adjoint group) in \cite[Table 1]{springer-tokyo} and using \cite{ccc} only for the few missing cases.}
\end{remark}

\section{Acknowledgements}
The authors are indebted to Kei Yuen Chan for communicating the statement of Lemma \ref{chan} and to Axel Hultman for pointing out the content of Remark~\ref{hultman}. This research was partially supported by Grants CPDA105885 and CPDA125818/12 of the University of Padova.

\end{document}